\documentclass[11pt]{amsart}

\theoremstyle{plain}
\newtheorem{thm}{Theorem}[section]

\newtheorem{cor}[thm]{Corollary}

\newtheorem{lemma}[thm]{Lemma}
\newtheorem{prop}[thm]{Proposition}

\theoremstyle{definition}
\newtheorem{defn}[thm]{Definition}

\newtheorem{example}[thm]{Example}

\newtheorem{genericitycond}[thm]{Genericity Condition}
\newtheorem{genericityconds}[thm]{Genericity Conditions}
\newtheorem{notation}[thm]{Notation}

\theoremstyle{remark}
\newtheorem{remark}[thm]{Remark}

\DeclareMathOperator{\mfg}{\mathfrak{g}}
\DeclareMathOperator{\mfh}{\mathfrak{h}}

\DeclareMathOperator{\jp}{\mathit{j}^{\prime}}
\DeclareMathOperator{\son}{\mathfrak{so}(n)}
\DeclareMathOperator{\sun}{\mathfrak{su}(n)}
\DeclareMathOperator{\sunp3}{\mathfrak{su}(n+3)}
\DeclareMathOperator{\sonp4}{\mathfrak{so}(n+4)}
\DeclareMathOperator{\spnp2}{\mathfrak{sp}(n+2)}
\DeclareMathOperator{\spn}{\mathfrak{sp}(n)}

\DeclareMathOperator{\Sunp3}{SU(n+3)}

\DeclareMathOperator{\ioj}{I_0(\it{g}_j)}
\DeclareMathOperator{\ioej}{I_0^e(\it{g_j})}

\DeclareMathOperator{\ioejp}{I_0^e(\it{g_{j^{\prime}}})}

\DeclareMathOperator{\Gnp}{\it{G}_{n+p}}
\DeclareMathOperator{\gn}{\mathfrak{g}_n}

\begin{document}

\newcommand{\spacing}[1]{\renewcommand{\baselinestretch}{#1}\large\normalsize}
\spacing{1}

\title{Isospectral Metrics and Potentials on Classical Compact Simple Lie Groups}

\author{Emily Proctor}

\address{Department of Mathematics and Statistics\\ Swarthmore College\\ Swarthmore, PA 19081}

\email{eprocto1@swarthmore.edu}

\thanks{Research partially supported by NSF grant DMS 0072534.}

\maketitle


\section{Introduction}\label{intro}

Given a compact Riemannian manifold $(M,g)$, the eigenvalues of the Laplace operator $\Delta$ form a discrete sequence known as the spectrum of $(M,g)$.  (In the case the $M$ has boundary, we stipulate either Dirichlet or Neumann boundary conditions.)  We say that two Riemannian manifolds are isospectral if they have the same spectrum.  For a fixed manifold $M$, an isospectral deformation of a metric $g_0$ on $M$ is a continuous family $\mathcal{F}$ of metrics on $M$ containing $g_0$ such that each metric $g\in\mathcal{F}$ is isospectral to $g_0$.  We say that the deformation is nontrivial if none of the other metrics in $\mathcal{F}$ are isometric to $g_0$ and that the deformation is multidimensional if $\mathcal{F}$ can be parameterized by more that one variable.  For two functions $\phi,\psi\in C^{\infty}(M)$, we say that $\phi$ and $\psi$ are isospectral potentials on $(M,g)$ if the eigenvalue spectra of the Schr\"odinger operators $\hbar\Delta +\phi$ and $\hbar\Delta + \psi$ are equal for any choice of Planck's constant $\hbar$. 

In this paper, we prove the existence of multiparameter isospectral deformations of metrics on $SO(n)$ $(n = 9$ or $n\geq 11)$, $SU(n)$ $(n\geq 8)$, and $Sp(n)$ $(n\geq 4)$.  For these examples, we follow a metric construction developed by Schueth who had given one-parameter families of isospectral metrics on orthogonal and unitary groups.  Our multiparameter families are obtained by a new proof of nontriviality establishing a generic condition for nonisometry of metrics arising from the construction.  We also show the existence of non-congruent pairs of isospectral potentials and nonisometric pairs of isospectral conformally equivalent metrics on $Sp(n)$ for $n\geq 6$.  

The industry of producing isospectral manifolds began in 1964 with Milnor's pair of 16-dimensional isospectral, nonisometric tori \cite{Mi}.  Several years later, in the early
1980's, new examples began to appear sporadically (e.g. \cite{Vi}, \cite{Ik1}, \cite{GW2}).   These isospectral constructions were ad hoc and did not appear to be related until 1985, when Sunada began developing the first unified approach for producing isospectral manifolds.  The method described a program for taking quotients of a
given manifold so that the resulting manifolds were isospectral.  Sunada's original theorem and subsequent generalizations (\cite{Bd1}, \cite{Bd2}, \cite{DG2}, \cite{Pe1},
\cite{Sut}) explained most of the previously known isospectral
examples and led to a wide variety of new examples.  See, for
example, \cite{BGG}, \cite{Bus}, and \cite{GWW}.

In 1993 Gordon produced the first examples of closed isospectral manifolds with different local geometry~\cite{Go93} and then, in a series of papers, generalized the construction to the following principle based on torus actions.

\begin{thm}\label{submersionthm} Let $T$ be a torus and suppose $(M,g)$ and
  $(M^{\prime},g^{\prime})$ are two principal $T$-bundles such that
  the fibers are totally geodesic flat tori.  Suppose that for any subtorus $K\subset T$ of
  codimension 0 or 1, the quotient manifolds $(M/K, \overline{g})$ and
  $(M^{\prime}/K, \overline{g}^{\prime})$, where $\overline{g}$ and
  $\overline{g}^{\prime}$ are the induced submersion metrics, are isospectral.
  Then $(M,g)$ and $(M^{\prime},g^{\prime})$ are isospectral.
\end{thm}

Gordon's initial application of
Theorem~\ref{submersionthm} was to give a sufficient condition for two
compact nilmanifolds (discrete quotients of nilpotent Lie groups) to
be isospectral.  In 1997, Gordon and Wilson furthered the development of the submersion
technique when they constructed the first examples of
\begin{it}continuous\end{it} families of isospectral manifolds with
different local geometry \cite{GW1}.  The base manifolds were products of
$n$-dimensional balls with $r$-dimensional tori ($n\geq 5$, $r\geq 2$),
realized as domains within nilmanifolds.  Here Gordon and Wilson
proved a general principle for local nonisometry within their construction.  They were also able to
exhibit specific examples of isospectral deformations of manifolds
with boundary for which the eigenvalues of the Ricci
tensor (which, in this setting, were constant functions on each
manifold) deformed nontrivially. 
It was later proven in \cite{GGSWW}
that the boundaries $S^{n-1}\times T^r$ of the manifolds in \cite{GW1} were
also examples of isospectral manifolds.  These were closed manifolds
which were not locally homogeneous.  A general abstract principle was
given for nonisometry but specific examples were also produced for
which the maximum
scalar curvature changed throughout the deformation, thereby proving
maximal scalar curvature is not a spectral invariant.

Expanding on the ideas of \cite{GGSWW}, Schueth produced the first
examples of simply connected closed isospectral manifolds.  In fact,
Schueth even produced continuous families of such manifolds
\cite{Sch99}.  Schueth's basic principle was to embed the torus $T^2$
into a larger, simply connected Lie group (e.g. $SU(2)\times
SU(2)\simeq S^3\times S^3$) and extend the metric in order to find
isospectral metrics on $S^4\times SU(2)\times SU(2)\simeq S^4\times
S^3\times S^3$.  Since the torus was embedded in the group, the torus action on
the manifold was the natural group action.  Schueth's examples were not
locally homogeneous.  For these examples, the critical values of the
scalar curvature changed throughout the deformation, proving the
manifolds were not locally isometric.  Furthermore, by examining heat
invariants related to the Laplacian on one-forms, Schueth was able to
prove that these examples were isospectral on functions but not on
one-forms.  

Schueth continued to capitalize on the notion of embedding the torus
in a larger group in \cite{Sch}.  In this
case, Schueth specialized Gordon's Theorem~\ref{submersionthm} to
compact Lie groups (see Theorem~\ref{dorotheeisospec}) in order to produce one-dimensional isospectral
deformations of each of $SO(n)\times T^2$ ($n\geq 5$), Spin$(n)\times T^2$, ($n\geq 5)$, $SU(n)\times
T^2$ ($n\geq 3$), $SO(n)$ ($n\geq
8$), Spin$(n)$ ($n\geq 8$), and $SU(n)$ ($n\geq 6$).  Here the metrics
were left-invariant so the manifolds were homogeneous.  As with many
previous examples, Schueth's metrics were constructed from linear maps
$j$ into the Lie algebra of the Lie group in question.  In order to
prove nonisometry, Schueth expressed the norm of the Ricci tensor in
terms of the associated linear map and chose her linear maps so that the norm of the Ricci tensor varied through the deformation. 

These particular examples of Schueth's were the inspiration for the first part of this
paper.  We will use Schueth's specialization of
Theorem~\ref{submersionthm} to produce our metrics.  However, in this
paper we will produce \begin{it}multidimensional\end{it} families of
metrics and will develop a general nonisometry principle for families of metrics
arising from linear maps according to Schueth's construction.  Furthermore, we will expand the class
of Lie groups for which such exist to include all of the classical
compact simple Lie groups of large enough dimension.

More recently, Gordon and Schueth have constructed conformally
equivalent metrics $\phi_1g$ and $\phi_2g$ on spheres $S^n$ and balls $B^{n+1}$ ($n\geq 7$)
and on $SO(n)$ ($n\geq 14$), Spin$(n)$ ($n\geq 14$), and $SU(n)$ ($n\geq 9$) \cite{GSch02}.  They also showed that the
conformal factors $\phi_1$ and $\phi_2$ were isospectral potentials for
the Schr\"{o}dinger operator $\hbar^2\Delta + \phi$ on each of these manifolds. 
In this paper, we extend their result to include $Sp(n)$ for $(n\geq 6)$.

The outline of the paper is as follows.  In Section~\ref{background} we describe the metrics and potentials which we will use and state the theorems by Gordon and Schueth which we will apply.  Next, in Section~\ref{examples} we give our examples of multiparameter isospectral deformations of metrics on the classical compact simple Lie groups.  Section~\ref{nonisometry} is devoted to proving the nontriviality of these examples.  Finally, in Section~\ref{spn} we give our examples of non-congruent isospectral potentials and nonisometric conformally equivalent isospectral metrics on $Sp(n)$ for $n\geq 6$.

The author is pleased to thank Carolyn Gordon and Dorothee Schueth for many helpful conversations about this work.


\vspace{6mm}

\section{Metric and Potential Constructions}\label{background}

In this section, we describe the metrics and potentials that we will consider for the remainder of the paper.  The constructions, which are due to Schueth and Gordon \cite{Sch},\cite{GSch02}, are based on linear maps.

Consider a Lie group $G$ with Lie algebra $\mfg$ and bi-invariant metric $g_0$.   By \begin{it}torus\end{it}, we mean a nontrivial, compact, connected abelian Lie group.  Suppose that $H<G$ is a torus with Lie algebra $\mfh$ and that $K<G$ is a closed connected subgroup with Lie algebra $\mathfrak{k}$.  Assume $\mfh$ is $g_0$-orthogonal to $\mathfrak{k}$ and $[\mfh,\mathfrak{k}]=0$.

\begin{notation}\label{metricnotation}Given a linear map $j:\mfh\to\mathfrak{k}\subset \mfg$, we define $j^t:\mfg\to \mfh$ by $g_0(j^t(X), Z) = g_0(X,j(Z))$
for all $X\in\mfg$, $Z\in\mfh$.  In other words, $j^t$ is the
$g_0$-transpose of $j$.  We then have an inner product $g_j$ on $\mfg$
given by $g_j = (Id + j^t)^*g_0$.   Let $g_j$ also denote the left-invariant metric on $G$ that is associated to this inner product.
\end{notation}

Notice that $g_j$ differs from $g_0$ only on
$\mathfrak{k}\oplus\mfh$, where we have used the linear map $j$ to redefine
orthogonality.  In particular, $j$ determines a subspace $S=\{X-j^t(X)|X\in\mfg\}$ which is
$g_j$-orthogonal to $\mfh$ and such that $g_j$ restricted to $S$ is
linearly isometric to $g_0$ restricted to $\mathfrak{k}$ via the map
$X-j^t(X)\mapsto X$.

Recall that a Lie algebra is compact if it is the Lie algebra of a compact Lie group.

\begin{defn}\label{isospecequivdefn}
Let $\mfg$ be a compact Lie algebra with associated Lie group $G$.  Let $\mfh$ be a real inner product space.  Suppose $j,\jp:\mfh\to\mfg$ are linear maps.  We say that $j$ and $\jp$ are \begin{it}isospectral\end{it}, denoted $j\sim \jp$, if for each $z\in\mfh$, there exists $A_z\in G$ such that $j(z) = \mathrm{Ad}(A_z)\jp(z)$.  We say that $j$ and $\jp$ are \begin{it}equivalent\end{it}, denoted $j\simeq \jp$, if there exists $A\in G$ and $C\in O(\mfh)$ such that $j(z) = \mathrm{Ad}(A)\jp(C(z))$ for all $z\in\mfh$.
\end{defn}

\begin{remark}\label{matrixconjugation}Note that in the case $G= SO(n)$, $SU(n)$, or
  $Sp(n)$, the map $\mathrm{Ad}(A):\mfg\to\mfg$ is given by matrix
  conjugation.  Thus we may rewrite the isospectrality condition as $A_zj(z)A_z^{-1} =\jp (z)$ and the equivalence condition as $Aj(z)A^{-1} = \jp (C(z))$.
\end{remark}

\begin{remark}\label{altdefinition}We use the definition of equivalence which was introduced in \cite{GW1} for the case $\mfg=\son$ and in \cite{Sch} for the case $\mfg=\sun$.  This differs slightly from the definition cited in \cite{GSch02} by Gordon and Schueth.   Gordon and Schueth's definition states that $j$ and $\jp$ are equivalent if there exists $C\in O(\mfh)$ and \begin{it}any\end{it} automorphism $\phi$ of $\mfg$ such that $j(z)=\phi( \jp(C(z)))$ for all $z\in\mfh$.  This means that Gordon and Schueth's definition is less restrictive except in the cases of $\son$ ($n$ odd) and $\spn$ where every automorphism of $\mfg$ is an inner automorphism by some element of $SO(n)$ or $Sp(n)$ respectively.  
\end{remark}

The following theorem by Schueth is a specialization of Gordon's
submersion theorem (Theorem~\ref{submersionthm}).

\begin{thm}\label{dorotheeisospec} \textnormal{(\cite{Sch})}  Let $G$ be a compact Lie group with Lie algebra $\mfg$, and let
$g_0$ be a bi-invariant metric on $G$.  Let $H< G$ be a torus in
$G$ with Lie algebra $\mfh\subset \mfg$.  Denote by $\mathfrak{u}$ the
$g_0$-orthogonal complement of the centralizer $\mathfrak{z}(\mfh)$ of
$\mfh$ in $\mfg$.  Let $\lambda$, $\lambda^{\prime}:\mfg\to\mfh$ be
two linear maps with $\lambda_{|_{\mfh\oplus\mathfrak{u}}}
=\lambda^{\prime}_{|_{\mfh\oplus\mathfrak{u}}} = 0$ which satisfy: For
every $z\in \mfh$ there exists $A_z\in G$ such that $A_z$ commutes
with $H$ and $\lambda^{\prime}_z = \mathrm{Ad}(A_z)^*\lambda_z$, where
$\lambda_z:=g_0(\lambda(\cdot),z)$ and
$\lambda^{\prime}_z:=g_0(\lambda^{\prime}(\cdot),z)$.  Denote by
$g_{\lambda}$ and $g_{\lambda^{\prime}}$ the left-invariant metrics on
$G$ which correspond to the scalar products $(\mathrm{Id}
+\lambda)^*g_0$ and $(\mathrm{Id} + \lambda^{\prime})^*g_0$ on
$\mfg$.  Then $(G,g_{\lambda})$ and $(G,g_{\lambda^{\prime}})$ are isospectral.
\end{thm}

In particular, if $j,\jp:\mfh\to\mathfrak{k}\subset\mfg$ are isospectral maps, then letting $\lambda=j^t$
and $\lambda^{\prime}=\jp^t$, we may conclude
that the metrics $g_j$ and
$g_{\jp}$ on $G$ described above are isospectral.  

We have a similar theorem for producing pairs of isospectral potentials and pairs of conformally equivalent isospectral metrics.

\begin{thm}\label{CandDtheorem}\textnormal{(\cite{GSch02})}\label{potentialconformal}  Let $G$ be a compact Lie group with Lie algebra $\mfg$, let $P$ be a compact Lie subgroup with Lie algebra of the form $\mathfrak{p} = \mathfrak{k}\oplus\mathfrak{k}$ for some Lie algebra $\mathfrak{k}$, and let $H<G$ be a torus with Lie algebra $\mfh$.  Suppose that $[\mathfrak{p},\mfh]=0$ and that $\mfh$ is orthogonal to $\mathfrak{p}$ with respect to a bi-invariant metric $g_0$ on $G$.  Let $j_1, j_2:\mfh\to\mathfrak{k}$ be isospectral linear maps as in Definition~\ref{isospecequivdefn}.  Define $j:\mfh\to\mathfrak{k}\oplus\mathfrak{k}=\mathfrak{p}$ by $j(Z)=(j_1(Z),j_2(Z))$.  Denote by $g_j$ the associated left-invariant metric on $G$.  Let $\phi$ be a smooth function on $G$ which is right invariant under $H$ and invariant under conjugation by elements of $P$.  Suppose that there exists an isometric automorphism $\tau$ of $(G,g_0)$ such that $\tau_{|_{H}}=\mathrm{Id}$ and such that $\tau_*$ restricts to the map $(X,Y)\mapsto (Y,X)$ on $\mathfrak{k}\oplus\mathfrak{k} = \mathfrak{p}\subset\mfg$.  Then:
\begin{enumerate}
\item $\phi$ and $\tau^*\phi$ are isospectral potentials on $(G,g_j)$.
\item If, in addition, $\phi$ is positive then $\phi g_j$ and $(\tau^*\phi)g_j$ are conformally equivalent isospectral metrics on $G$.
\end{enumerate}
\end{thm}


\vspace{6mm}

\section{Examples of Isospectral Deformations of Metrics on Lie Groups}\label{examples}

We now apply the material in Section~\ref{background} to produce examples of isospectral deformations of metrics on Lie groups.  All of our examples arise from the following theorem.  In Section~\ref{nonisometry} we prove the nontriviality of the examples.  

\begin{thm}\label{families}
Suppose $\mfg$ is one of $\son\ (n=5, n\geq 7)$ , $\sun\ (n\geq4)$, or
$\spn\ (n\geq 3)$ with associated group $SO(n)$, $SU(n)$, or $Sp(n)$ respectively.  Suppose $\mfh$ is the Lie algebra of the two-dimensional torus.  Let $L$ be the space of all linear maps
$j:\mfh\rightarrow\mfg$.   Then there exists a Zariski open set $\mathcal{O}
\subset L$ such that each $j_0 \in \mathcal{O}$ is contained in a
continuous $d$-parameter family of linear maps which are isospectral but pairwise not
equivalent.  Here $d$ depends on $\mfg$ as follows:


\begin{table}[h]
\begin{center}
\begin{tabular}{ | l | l | }
\hline
\begin{bf}$\mfg$\end{bf} & \begin{bf}$d$\end{bf} \\ \hline
$\son$ & $d\geq n(n-1)/2 - [\frac{n}{2}]([\frac{n}{2}] + 2)$ \\ \hline
$\sun$ & $d \geq n^2 -1 - \frac{n^2 + 3n}{2}$ \\ \hline
$\spn$ & $d\geq n^2 - n$ \\ \hline
\end{tabular}
\end{center}
\end{table}
(The symbol $[\frac{n}{2}]$ denotes the largest integer less than or equal to $\frac{n}{2}$.)
Note that for $\son$, $d>1$ when $n=5$ or $n\geq 7$.  For $\sun$,
$d=1$ when $n=4$ and $d>1$ when $n\geq 5$.  For $\spn$, $d>1$ when
$n\geq 3$.
\end{thm}

\begin{remark}\label{CandEproof} This theorem was originally proven by Gordon and Wilson in \cite{GW1} for the case of $\son$ with the associated group $O(n)$.   However for fixed $j_0\in\mathcal{O}$, since the $d$-parameter family of linear maps which are isospectral to $j_0$ is continuous, we have that for each $z\in\mfh$, the family $\{j(z)\ |\ j\ \textrm{is an element of the \textit{d}-parameter family}\}$ is the orbit of $j_0(z)$ under the adjoint action of a continuous set of elements of $O(n)$.  The identity is contained in this set so therefore the set is in fact contained in $SO(n)$.  Thus we have that each of Gordon and Wilson's families consists of maps which are isospectral via $SO(n)$.  We will use this in our Example~\ref{multison}.  On the other hand, for any pair of maps within one of Gordon and Wilson's families, there is no element of $O(n)$ which makes them equivalent.  This is stronger than pairwise nonequivalence via $SO(n)$.  Indeed, for $n=5,6,7$ and $n\geq 9$, the automorphism group of $\son$ is contained in $\{Ad(A)\ |\ A\in O(n)\}$.  Thus we see that except for $n=8$, Gordon and Wilson's families consist of linear maps which are not equivalent even according to the definition given in \cite{GSch02} (cf. Remark~\ref{altdefinition}).  This will prove useful at the end of Section~\ref{nonisometry} when we prove nontriviality of our isospectral examples. 
\end{remark}

\begin{remark}\label{so6} Though Gordon and Wilson's original proof for $\son$ broke down in the case $n=6$ (since $n(n-1)/2 - [\frac{n}{2}]([\frac{n}{2}] + 2) = 0$ for $n=6$), they explicitly exhibited one-parameter families of isospectral, nonequivalent linear maps $j:\mfh\to\mathfrak{so}(6)$.  
\end{remark}

\begin{remark}
In \cite{Sch}, Schueth gives examples of one-dimensional families of isospectral, nonisometric metrics on $SO(n)$ $(n\geq9)$, $\mathrm{Spin}(n)$ $(n\geq 9)$, and $SU(n)$ $(n\geq6)$ which arise from continuous one-dimensional families of isospectral linear maps from $\mfh$ to $\son$ ($n\geq 5$) and $\sun$ ($n\geq 3$), respectively.  However, her proof of nonisometry differs from ours. Schueth's proof has the advantage of showing nonisometry in a geometric way by using the norm of the Ricci tensor, while ours has the advantage of giving a generic condition for linear maps to produce nonisometric metrics.
\end{remark}

\begin{remark}\label{so5sp2}
Gordon and Wilson's proof extends in a straightforward way to $\sun$ and $\spn$, making obvious adjustments depending on the Lie algebra.  For full details see \cite{Pr}.   Here, for $\spn$, we proved the case where $n\geq 3$.  However, since $\mathfrak{sp}(2)$ is isomorphic to $\mathfrak{so}(5)$, the result is also true for $n=2$.  
\end{remark}

\medskip

 Let $T^r$ denote an $r$-dimensional torus.  

\begin{example}\label{multison} Isospectral deformations of metrics on $SO(n)$.

Notice that $SO(n)\times T^r$ is contained as a subgroup of $SO(n)\times SO(2r)$ which is itself contained as a subgroup of $SO(n+2r)$ in the form of diagonal block matrices.  From Theorem~\ref{families} and Remark~\ref{CandEproof}, we have examples of one-dimensional families of pairwise isospectral linear maps $j:\mfh\to\mathfrak{so}(6)$ and multidimensional families of pairwise isospectral linear maps $j:\mfh\to\son$ for $n=5$ or $n\geq 7$.  Thus, according to the construction described in Section~\ref{background}, this gives one-dimensional families of isospectral metrics on $SO(10)$ and multidimensional families for $SO(9)$ and $SO(n)$ ($n\geq 11$).  By Theorem~\ref{dorotheeisospec}, within each family, the metrics are pairwise isospectral. 
\end{example}

\begin{example}\label{multispin} Isospectral deformations of metrics on Spin$(n)$.

By lifting from $SO(n)$ to Spin$(n)$, we may consider Spin$(n)\times T^r$ a subgroup of Spin$(n+2r)$.  Since the orbits of $\mathrm{Ad}(\mathrm{Spin}(n))$ in $\son$ are
  equal to the orbits of $\mathrm{Ad}(SO(n))$, if $j$ and $\jp$ are
  isospectral with respect to $SO(n)$, they are also isospectral with
  respect to Spin$(n)$.  Fix $j_0\in\mathcal{O}$.  By an argument similar to the one for
  $SO(n)$, we conclude that we have one-dimensional families of isospectral metrics on $\mathrm{Spin}(10)$ and multidimensional families of isospectral metrics on $\mathrm{Spin}(9)$ and $\mathrm{Spin}(n)$ for $n\geq 11$.\end{example}

\begin{example}\label{multisun} Isospectral deformations of metrics on $SU(n)$.

We have $SU(n)\times T^r$ contained as a subgroup of $SU(n)\times SU(r+1)$ which is contained in turn as a subgroup of $SU(n+r+1)$ .  From Theorem~\ref{families}, we have one-dimensional families of pairwise isospectral linear maps $j:\mfh\to\mathfrak{su}(4)$ and multidimensional families of pairwise isospectral linear maps $j:\mfh\to\sun$ for $n\geq 5$.  Thus we have one-dimensional isospectral deformations of metrics on $SU(7)$ and multidimensional families of isospectral metrics on $SU(n)$, $n\geq 8$.  
\end{example}

\begin{example}\label{multispn} Isospectral deformations of metrics on $Sp(n)$.

$Sp(n)\times T^r$ is contained as a subgroup of $Sp(n)\times Sp(r)$ which is contained as a subgroup of $Sp(n+r)$.   Theorem~\ref{families} and Remark~\ref{so5sp2} gives us multidimensional families of pairwise isospectral linear maps $j:\mfh\to\spn$ for $n\geq 2$.   Thus we have multidimensional families of isospectral metrics on $Sp(n)$ for $n\geq 4$.  
\end{example}


\vspace{6mm}

\section{Nonisometry of examples}\label{nonisometry}
In this section we prove nontriviality of Examples~\ref{multison} through \ref{multispn}.  Here we let $G_n$ denote one of $SO(n)$, Spin$(n)$, $SU(n)$, or $Sp(n)$ and $\mfg_n$ denote the associated Lie algebra $\son$, $\sun$, or $\spn$.  Recall that for Examples~\ref{multison} through \ref{multispn} we embedded $G_n\times T^2$ into a higher dimensional group.  For this section, we will refer to the higher dimensional group as $\Gnp$, where $p=4$ in the cases of $SO(n)$ and $\mathrm{Spin}(n)$, $p=3$ in the case of $SU(n)$, and $p=2$ in the case of $Sp(n)$.  Let $\ioj$ denote the identity component of the isometry group of
$(\Gnp, g_j)$ and let $\ioej$ denote the isotropy subgroup at $e$ of
$\ioj$.  For $x\in \Gnp$, denote left (resp. right) translation by $x$ by
$L_x$ (resp. $R_x$). 

\begin{thm}\textnormal{\cite{OT}}\label{Ochiai} Let $G$ be a compact, connected simple Lie group
  and $ds^2$ be a left-invariant Riemannian metric on $G$.  Then for each
  isometry $f$ contained in the
  identity component of the group of isometries of $(G,ds^2)$
  there exist $x,y\in G$ such that $f=L_x\circ R_y$. 
\end{thm}

In particular, for $\alpha\in\ioej$, we have that there exists some
$x\in\Gnp$ such that $\alpha$ is equal to conjugation of $\Gnp$ by
$x$.  Since $\alpha$ fixes the identity, at the Lie algebra level we
have that $\alpha_*$ is equal to $\mathrm{Ad}(x)$.

\begin{prop}\label{simplegroupisom} Suppose $G$ is a compact simple
  group with left-invariant metrics $g$ and $g^{\prime}$ neither of which is bi-invariant.  If $\mu:(G,g) \to (G,g^{\prime})$ is
  an isometry such that $\mu(e)=e$, then $\mu$ is an automorphism of $G$.
\end{prop}

\begin{proof}Since $\mu$ is an isometry, we have that $\ioej$ is isomorphic to $\ioejp$ via conjugation by $\mu$.  $G$ is compact so
  the isometry groups of $(G,g)$ and $(G,g^{\prime})$ are also
  compact.  Thus we may write the isometry
  group of $(G,g)$ as $G_1 \times G_2
  \times \dots \times G_s \times T/Z$ and the isometry group of
  $(G,g^{\prime})$ as $G_1^{\prime} \times
  G_2^{\prime} \times \dots \times G_t^{\prime}\times T^{\prime}/Z^{\prime}$
  where each $G_i^{(\prime)}$ is simple, $T^{(\prime)}$ is a torus, and $Z^{(\prime)}$ is
  central.  Each isometry group contains a copy of $G$ in the form of
  left translations.  Furthermore, since neither $g$ nor $g^{\prime}$
  is bi-invariant, each isometry group contains exactly one copy of
  $G$ by Theorem~\ref{Ochiai}.  Any isomorphism from the isometry group of $(G,g)$ to the
  isometry group of $(G,g^{\prime})$ must carry simple
  factors to simple factors.  Since $G$ is the only simple factor
  of its dimension, we have that any isomorphism carries $G$ to
  $G$.  This means that for every $x\in
  G$ there exists $x^{\prime}\in G$ such that $\mu
  L_x\mu^{-1} = L_{x^{\prime}}$.  Thus for each $x, y\in G$,
\begin{multline}
\mu(xy) = \mu(L_{xy})(e) = \mu(L_xL_y)(e) = \mu L_x \mu^{-1}\mu
L_y(e)\\ = L_{x^{\prime}}\mu L_y(e) = L_{x^{\prime}}L_{y^{\prime}}\mu (e) =
L_{x^{\prime}}L_{y^{\prime}}(e) = x^{\prime}y^{\prime} = \mu (x)\mu (y).
\end{multline} 
\end{proof}

\begin{remark}
If $\mu:(G,g)\to (G,g^{\prime})$ is an isometry which does not carry the identity to itself, then composing $\mu$ with $L_{\mu(e)^{-1}}$, we have an automorphism of $G$.
\end{remark}

The following corollary is immediate.

\begin{cor}\label{sunp3isom} Let $\Gnp$ be one of $SO(n+4)$, Spin$(n+4)$, $SU(n+3)$, or $Sp(n+2)$.  Given two nonzero linear maps, $j$ and $\jp$, suppose there exists
  an isometry $\mu : (\Gnp, g_j) \to (\Gnp, g_{\jp})$, where $g_j$ and $g_{\jp}$ are as in~\ref{metricnotation}.  Then
 $L_{\mu(e)^{-1}}\circ\mu$ is an
  automorphism of $\Gnp$.
\end{cor}

\begin{lemma}\label{t2equivalence} Let $\Gnp$ be one of $SO(n+4)$, Spin$(n+4)$, $SU(n+3)$, or $Sp(n+2)$, where $T^2$ is embedded, as earlier, in $G_p\subset G_{n+p}$.  Let $j,\jp:\mfh\to\mfg_n$ be nonzero linear maps and let $g_j$ and $g_{\jp}$ be as in~\ref{metricnotation}.  Suppose there exists an isometry
$\mu : (\Gnp, g_j) \to (\Gnp, g_{\jp})$ such that $\mu(e) = e$ and $\mu(T^2) =
T^2$.  Then there is an element $C\in O(\mfh)$ such that $j(z) =
\mu^{-1}_*\jp(Cz)$ for all $z\in\mfh$.
\end{lemma}

\begin{proof}
From Corollary~\ref{sunp3isom} we know that $\mu$ is an automorphism of $G_{n+p}$.  Thus $\mu_*$ maps left invariant vector fields to left invariant vector fields, i.e. $\mu_*:\mfg_{n+p}\to\mfg_{n+p}$.  If $\mu$ maps $T^2$ to itself, then it must isometrically map the Lie
algebra $\mfh$ to itself.  This implies that there is an
element $C\in O(\mfh)$ such that $\mu_*$ restricted to $\mfh$ is equal
to $C$.

Furthermore, since $\mu$ is an automorphism of $G_{n+p}$, if $\mu$ maps $T^2$
to itself in $\Gnp$, it must also isomorphically map the identity
component of the centralizer
of $T^2$ in $\Gnp$ to itself.  At the Lie algebra
level, direct calculation shows that the centralizer of $\mfh$ in
 
\begin{itemize}
\item $\sonp4$ is $\son\oplus\mfh$.  Thus the identity
  component of the centralizer of
  $T^2$ in $SO(n+4)$ is $SO(n)\times T^2$ and the identity component of the centralizer
  of $T^2$ in $\mathrm{Spin}(n+4)$ is
  $\mathrm{Spin}(n)\times T^2$. 
\item $\sunp3$ is
$\sun\oplus tu\oplus\mfh$ where
\begin{equation}\label{u}
u=
\left[
\begin{smallmatrix}
i/n & {} & {} & {} & {} & {}\\
{} & \ddots & {} & {} & {} & {}\\
{} & {} & i/n & {} & {} & {}\\
{} & {} & {} & -i/3 & {} & {}\\
{} & {} & {} & {} & -i/3 & {}\\
{} & {} & {} & {} & {} & -i/3
\end{smallmatrix}
\right] \in \sunp3.
\end{equation}  
Letting $U$ denote the
one-parameter subgroup associated to $u$, we have the identity
component of the centralizer of $T^2$ in $\Sunp3$ is $SU(n)\times U \times T^2$.
\item $\spnp2$ is $\spn\oplus\mfh$.  Thus the identity component of
  the centralizer of
  $T^2$ in $Sp(n+2)$ is $Sp(n)\times T^2$.
\end{itemize}

In each case, the identity component of the centralizer of $T^2$ is
the product of a simple group, $G_n$, with a torus.
Therefore $\mu(G_n) = G_n$ and $\mu_*$ is a Lie algebra automorphism
of $\mfg_n$.  Hence, for any $X\in\gn$, we have that $X -
j^t(X)\in\mfh^{\perp_{g_j}}$ is mapped to $\mu_*X - Cj^t(X)$.

On the other hand, since
$\mu_*X\in\gn$ and $Cj^t (X)\in\mfh$ and since
$\mu_*:\mfh^{\perp_{g_j}}\to\mfh^{\perp_{g_{\jp}}}$, it must be the case that $Cj^t (X)
= \jp^t (\mu_*X)$ for all $X\in\gn$.  Otherwise, $Cj^t
(X) = \jp^t(\mu_*X) +Z$ for some nonzero $Z\in\mfh$ depending
on $X$.  But in this case,
$\mu_*X-Cj^t (X) = \mu_*X - \jp^t (\mu_*X)
- Z$ which is not in
$\mfh^{\perp_{g_{\jp}}}$.  

Finally, taking transposes, we see that the condition
$j^t (X) = C^{-1}\jp^t(\mu_*X)$ for all $X\in\gn$ implies $j(z) =
\mu^{-1}_*\jp(Cz)$ for all $z\in\mfh$.
\end{proof}

\begin{remark}\label{innerautomorphism} From the proof of Lemma~\ref{t2equivalence}, we see that $\mu$ restricted to $G_n$ is an automorphism.  Suppose $\mu$ restricted to $G_n$ is an inner
automorphism so $\mu_*$ restricted to $\mfg_n$ is equal to
  $\mathrm{Ad}(A)$ for some $A\in G_n$.  Then by the proof, $j(z) =
  \mathrm{Ad}(A^{-1})\jp(Cz)$ for all $z\in\mfh$.  In other words, $j$ and
  $\jp$ are equivalent.
\end{remark}

\begin{genericitycond}\label{genericitycondition}We say that $j:\mfh\to\mfg_n$
  is generic if there are only finitely many $A\in G_n$ such that
  $j(z)=\mathrm{Ad}(A)j(z)$ for all $z\in\mfh$.
\end{genericitycond}

From the proofs found in \cite{GW1} and \cite{Pr}, the linear maps $j$ used in Examples~\ref{multison} through \ref{multispn} are generic.

\begin{lemma}\label{maximaltorus}
Let $\Gnp$ be one of $SO(n+4)$, Spin$(n+4)$, $SU(n+3)$, or $Sp(n+2)$.  Let $j:\mfh\to\gn$ be generic and let $g_j$ be the associated metric on
$\Gnp$. For $\Gnp$ equal to $SO(n+4)$, $\mathrm{Spin}(n+4)$, or $Sp(n+2)$, let $D$ be the group
of isometries of $(\Gnp,g_j)$ generated by the set $\{L_x\circ
R_{x^{-1}}|x\in T^2\}$.  For $\Gnp$ equal to $SU(n+3)$, let $D$ be the
group of isometries of $(\Gnp,g_j)$ generated by the set $\{L_x\circ
R_{x^{-1}}|x\in U\times T^2\}$, where $U$ is as in the proof of Lemma~\ref{t2equivalence}.  Then $D$ is a maximal torus in $\ioej$.
\end{lemma}

\begin{proof}
It is straightforward to check that $D \subset \ioej$.  Recall that every element of $\ioej$ is of the form $L_x\circ R_{x^{-1}}$ for some
$x\in\Gnp$.  Let $C(\Gnp)$ denote the finite center of $\Gnp$.
We identify $\ioej$ with a subgroup of $\Gnp/C(\Gnp)$ via the
map which sends $L_x\circ R_{x^{-1}}$ to the coset of $x$ in
$\Gnp/C(\Gnp)$.  Under this correspondence, we consider $D$ a subgroup of $\Gnp/C(\Gnp)$. 

Suppose that $\{y_t\ \vert\ t\in (-\epsilon,\epsilon)\}$ is a continuous family of elements of $G_{n+p}$ with $y_0=e$.  Furthermore, suppose that for each $t$, $L_{y_t}\circ R_{y_t^{-1}}$ is an element of $\ioej$ which commutes with $D$.  If $L_{y_t}\circ R_{y_t^{-1}}$ commutes with $D$, then under the identification of $D$ with a subgroup of $G_{n+p}/C(\Gnp)$, for each $x\in T^2$ (resp. $U\times T^2$) $y_t^{-1}xy_t=xz_t$ for some $z_t\in C(\Gnp)$.  But $C(\Gnp)$ is discrete so it must be the case that $z_t=e$ for all $t$ and all $x\in T^2$.  This implies that for any $t$ when $y_t$ acts by isometry (i.e. conjugation) on $\Gnp$, it fixes $T^2$ pointwise.  Thus the continuous family $\{y_t\ \vert\ t\in(-\epsilon, \epsilon)\}$ is contained in the identity component of the centralizer of $T^2$ of $\Gnp$.

Since $y_t$ is contained in the identity component of the centralizer of $T^2$ for each $t$, we have that for $SO(n+4)$, $\mathrm{Spin}(n+4)$, and $Sp(n+2)$, $y_t$ is
equal to a product $A_tZ_t$ and for $SU(n+3)$, $y$ equals $A_tUZ_t$ for some $A_t\in G_n$, $Z_t\in
T^2$.  In this case, $(L_{y_t}\circ R_{y_t^{-1}})_*$ restricted to $\mfg_n$ equals $\mathrm{Ad}(A_t)$.   By the proof of Lemma~\ref{t2equivalence}, $j(z) = \mathrm{Ad}(A_t^{-1})j(Cz)$ for some $C\in O(\mfh)$.  But since $L_{y_t}\circ R_{y_t^{-1}}$ fixes $T^2$ pointwise, we have that $C$ is equal to the identity so
\begin{equation}\label{finite}
j(z) =\mathrm{Ad}(A_t^{-1})j(z)\hspace{2mm}\mathrm{for}\ \mathrm{all}\ z\in\mfh.
\end{equation}
But by the genericity of $j$, there are only finitely many $A_t$ for which Equation~\ref{finite} holds and thus only finitely many $A_t$ such that $y_t=A_tZ_t$ (resp. $A_tUZ_t$) for some $Z_t\in T^2$.  Since our family $\{y_t\ \vert\ t\in(-\epsilon,\epsilon)\}$ is continuous, it must be that $A_t$ is the identity for all $t\in (-\epsilon,\epsilon)$.  Therefore $\{y_t\ \vert\ t\in(-\epsilon,\epsilon)\}\subset T^2$.  In other words, $D$ is not contained in a higher dimensional connected torus and hence is a maximal torus in $\ioej$.  

\end{proof}

\begin{thm}\label{bigtheorem}Let $\Gnp$ be one of $SO(n+4)$, Spin$(n+4)$, $SU(n+3)$, or $Sp(n+2)$.  Let $j$ and $\jp$ be generic linear maps such 
 that $\mu: (\Gnp, g_j) \to (\Gnp, g_{\jp})$ be an isometry.
  Then there exists an element $C\in O(\mfh)$ such that $j(z) =
 \mu_*^{-1}\jp(Cz)$ for all $z\in\mfh$.  By
  Remark~\ref{innerautomorphism}, if $\mu$ restricts to an
  inner automorphism of $G_n$, then $j$ and $\jp$ are equivalent.
\end{thm}

\begin{proof}
Suppose that $\mu : (\Gnp,g_j)\to (\Gnp,g_{\jp})$ is an isometry.
  We may assume that $\mu(e) = e$.  By Lemma~\ref{t2equivalence}, it suffices to show that $\mu(T^2)=T^2$.

We know that $\ioej$ is isomorphic
to $\ioejp$ via conjugation by $\mu$.  According to Lemma~\ref{maximaltorus},
$D$ is a maximal torus in $\ioej$ and so the isomorphism carries $D$
to a maximal torus in $\ioejp$.  All maximal tori in a compact Lie
group are conjugate so, after possibly composing $\mu$ with an
element of $\ioejp$, we may assume that conjugation by $\mu$ carries $D$
to the similarly defined set in $\ioejp$.

\begin{itemize}
\item For $SO(n+4)$, Spin$(n+4)$, and $Sp(n+2)$, this implies that for any $a\in
  T^2$, $\mu\circ L_a \circ
R_{a^{-1}} \circ \mu^{-1}= L_b\circ R_{b^{-1}}$ for some $b\in T^2$.
On the other hand, by Corollary~\ref{sunp3isom}, we know $\mu$ is an
automorphism of $\Gnp$.  Thus, for any $x\in\Gnp$, 
\begin{equation*}
\mu\circ L_a \circ R_{a^{-1}}\circ \mu^{-1} (x) = \mu(a\mu^{-1}(x)a^{-1}) =
\mu (a)x\mu^{-1}(a) = L_{\mu(a)}\circ R_{\mu(a^{-1})}(x).
\end{equation*}
In other words, $\mu(a)=bz$ for some $z\in C(\Gnp)$.  For each of
$SO(n+4)$, Spin$(n+4)$, and $Sp(n+2)$, $C(\Gnp)$ is finite.  Since $\mu$ is continuous and
since $\mu(e)=e$, we have that $\mu(T^2)=T^2$.

\item For $SU(n+3)$, a similar argument implies that $\mu(U\times T^2) = U\times T^2$ and thus, at the Lie algebra
level, $\mu_*$ maps $tu \oplus \mfh$ to $tu\oplus \mfh$.  

Now consider
$\mu$ as an automorphism.  The automorphism group of $\Sunp3$ is
generated by the inner automorphisms and one outer automorphism,
namely complex conjugation.  At the Lie algebra level, conjugating an element $X\in\sunp3$ by any
element of $\Sunp3$ preserves the eigenvalues of $X$.  In particular,
$u$ has eigenvalue $i/n$ with multiplicity $n$ and eigenvalue $-i/3$
with multiplicity $3$.  No other element of $tu\oplus\mfh$ has the
same eigenvalues.  Therefore each inner automorphism of $\sunp3$ maps
$u$ to $u$.  Similarly, at the Lie algebra level, the outer automorphism of
$\Sunp3$ negates the eigenvalues of $X\in\sunp3$.  For each
$t\in\Bbb{R}$, this sends $tu$ to $-tu$.  Thus the vector space
spanned by $u$ is fixed.

Since $\mu$ is an isometry, $\mu_*$ maps the $g_j$-orthogonal
complement of the space spanned by $u$ to the $g_{\jp}$-orthogonal
complement of the space spanned by $u$. Thus $\mu_*(\mfh) = \mfh$ and therefore $\mu(T^2)=T^2$.  
\end{itemize}
\end{proof}

\begin{thm}\label{finitelymany}Let $\mfg_n$ be one of $\son$, $\sun$, or $\spn$.  Suppose $j_0:\mfh\to\mfg_n$ is
  contained in a family of generic linear maps which are pairwise nonequivalent.   Then there are at most finitely other maps $j$ contained in the family such that $g_{j_0}$ and $g_j$ are isometric.
\end{thm}
\begin{proof}

Suppose $j$ and $\jp$ are two generic linear maps such that $g_j$ and $g_{\jp}$ are
both isometric to $g_{j_0}$.  From Theorem~\ref{bigtheorem} we have
Lie algebra automorphisms $\alpha$, $\alpha^{\prime}$ of $\mfg_n$ and elements $C,
C^{\prime}$ of $O(\mfh)$ such that 
\begin{equation}
j(z) = \alpha j_0 (Cz)
\end{equation}
and 
\begin{equation}
\jp(z) = \alpha^{\prime} j_0 (C^{\prime}z)
\end{equation}
for all $z\in\mfh$.

If $\alpha$ and $\alpha^{\prime}$ are in the same coset of
$\mathrm{Aut}(\mfg)/\mathrm{Aut}^0(\mfg)$, then they differ by an
inner automorphism.  But in this case $j$ and $\jp$ are equivalent.
Since our family consists of pairwise nonequivalent linear maps, $j$
is equal to
$\jp$.  The theorem now follows from the fact that for each of $\son$,
 $\sun$, and $\spn$,
  $\mathrm{Aut}(\mfg)/\mathrm{Aut}^0(\mfg)$ is finite.
\end{proof}

Finally we may apply Theorem~\ref{finitelymany} to Examples~\ref{multison} through \ref{multispn}.   For Examples \ref{multison} and \ref{multispin}, let $j_0$ be a linear map in the set $\mathcal{O}$ from Theorem~\ref{families} and let $\mathcal{F}_{\son}$ be a $d$-parameter family of linear maps which are isospectral but pairwise nonequivalent.  Since the automorphism group of $\son$ ($n=5, 6, 7$, and
$n\geq 9$) is contained in $\{\mathrm{Ad}(A)|A\in O(n)\}$, we have
from the proof of Theorem~\ref{finitelymany} and Remark~\ref{CandEproof} that no two elements of
$\mathcal{F}_{\son}$ give rise to isometric metrics on $SO(n)$ or $\mathrm{Spin}(n)$.  Thus we have an isospectral deformation of $g_{j_0}$.  For $\mathfrak{so}(8)$, since
the cardinality of $\mathrm{Aut}(\mathfrak{so}(8))/\mathrm{Aut}^0(\mathfrak{so}(8))$ is three, for any element of
$\mathcal{F}_{\mathfrak{so}(8)}$ there are at most two other elements which could give rise
to isometric metrics on $SO(12)$ or $\mathrm{Spin}(12)$.  For fixed $j_0\in\mathcal{O}$, we may choose
$\mathcal{F}_{\mathfrak{so}(8)}$ small enough that no other element of $\mathcal{F}_{\mathfrak{so}(8)}$
produces a metric isometric to $g_{j_0}$, thereby obtaining an
isospectral deformation of $g_{j_0}$.

The analysis of Examples~\ref{multisun} and \ref{multispn} is similar.  Recall that the automorphism group of $\sun$ for $n\geq 4$ is generated by inner automorphisms and the outer automorphism which takes an element to its complex conjugate.  By choosing $j_0$ in $\mathcal{O}$, we have a $d$-dimensional isospectral deformation of $g_{j_0}$ such that for any metric other than $g_{j_0}$ within the deformation, there is at most one other isometric metric contained in the deformation.  Finally, the automorphism group of $\spn$ consists entirely of inner automorphisms.  Thus for $n\geq 4$, we have
multiparameter isospectral deformations of metrics on $Sp(n)$ such
that no two metrics in a given deformation are isometric.

Thus we have produced isospectral deformations of metrics on each of
$SO(n)$ ($n\geq9$), Spin$(n)$ ($n\geq 9$), $SU(n)$
($n\geq 7$), and $Sp(n)$ ($n\geq 4$).  Except for low dimensions, all of these deformations are multidimensional.


\vspace{6mm}

\section{Examples of Isospectral Potentials and Conformally Equivalent Isospectral Metrics on $Sp(n)$}\label{spn}

Here we apply Theorem~\ref{CandDtheorem} to $Sp(n)$.  For a particular choice of $\phi$ we will find that our pairs of isospectral potentials are non-congruent and our pairs of conformally equivalent isospectral metrics are nonisometric.  As above, we let $T^r$ denote an $r$-dimensional torus with Lie algebra $\mfh$.

\begin{example}
First, consider $T^r$ as an embedded subtorus of $Sp(r)$ and consider $Sp(n)\times Sp(n)\times Sp(r)$ as a block-diagonal subgroup of $Sp(2n+r)$.  (Note: We are thinking of elements of $Sp(2n+r)$ as complex matrices with $(2n+r)^2$ $2\times 2$ blocks.)  Suppose $\tau$ is the automorphism of $Sp(2n+r)$ given by conjugation by the matrix 
\begin{equation}
c=
\begin{bmatrix}
0 & \mathrm{Id}_{2n} & 0\\
\mathrm{Id}_{2n} & 0 & 0\\
0 & 0 & \mathrm{Id}_{2r}
\end{bmatrix}.
\end{equation}

From Theorem~\ref{families} we have continuous families of isospectral linear maps from the Lie algebra $\mfh$ of $T^2$ into $\spn$ for $n\geq 2$.  Suppose that $j_1$ and $j_2$ are two elements of a particular family satisfying the following genericity conditions.  

\begin{genericityconds}\label{GSchgenericity}.

a. The kernel of the map $j=(j_1,j_2):\mfh\to\spn\oplus\spn$ is trivial.

b. The image of $j$ has trivial centralizer in $\spn\oplus\spn$.
\end{genericityconds}

Notice that condition~\ref{GSchgenericity}a is a very mild condition which is easily satisfied.  Suppose that condition~\ref{GSchgenericity}b is not satisfied by the maps $j_1$ and $j_2$.  This would imply that there is at least a one-dimensional family of elements $X$ in $\spn$ such that $[X,j_1(Z)]=0$ for all $Z\in\mfh$.  But then it would follow that there exists a one-dimensional family of elements $A$ in $Sp(n)$ such that $\mathrm{Ad}(A)j_1(Z)=j_1(Z)$ for all $Z\in\mfh$.  Thus, if $j_1, j_2$ are generic in the sense of Condition~\ref{genericitycondition} then Condition~\ref{GSchgenericity}(b) is automatically satisfied.

Let $\phi$ be any smooth function on $G$ which is right invariant uner $T^2$ and invariant under conjugation by elements of $Sp(n)\times Sp(n)$.  Then, following the construction in Theorem~\ref{CandDtheorem}, we have isospectral potentials $\phi$ and $\tau^*\phi$ on $Sp(2n+2,g_j)$ for $n\geq 2$.  Furthermore, if $\phi$ is positive then $\phi g_j$ and $(\tau^*\phi)g_j$ are conformally equivalent isospectral metrics on $Sp(2n+2)$.

We now choose $\phi$ so that $\phi$ and $\tau^*\phi$ are not congruent on $Sp(2n+2, g_j)$ and so that $\phi g_j$ and $(\tau^*\phi)g_j$ are not isometric.  We follow Gordon and Schueth's construction and suppose that each matrix in $Sp(2n+r)$ is written
\begin{equation}
X = 
\begin{bmatrix}
A & B & C\\
D & E & F\\
H & J & L
\end{bmatrix}
\end{equation}
where $A, B, C, D$ are $2n\times 2n$ matrices, $L$ is a $2r\times 2r$ matrix, $C, F$ are $2n\times 2r$ matrices, and $H,J$ are $2r\times 2n$ matrices.
Let $c_1>c_2>0$ and define $\phi$ by 
\begin{equation}
\phi(X) :=\, \mathrm{exp}(c_1\ \mathrm{Re\ det}\ A +c_2\ \mathrm{Re\ det}\ E).
\end{equation}

This function $\phi$ is almost exactly the same as the function $\phi$ used for $SO(2m+2r)$ in \cite{GSch02} except that here we have taken the real parts of the determinants of $\det A$ and $\det E$ to account for the fact that matrices in $Sp(n)$ have complex entries.  In their proof, Gordon and Schueth had to treat the cases of $SO(2m+2r)$ and $SU(2m + r+ 1)$ separately because of the elements 
$\left[\begin{smallmatrix}
\alpha\mathrm{Id}_m & &\\
& \beta\mathrm{Id}_m &\\
& & \gamma\mathrm{Id}_{r+1}
\end{smallmatrix}\right]$
in $SU(2m+r+1)$ (cf. proof of Lemma~\ref{t2equivalence} above).  The analogous element does not exist in $Sp(2n+r)$.
Thus the proofs that $\phi$ and $\tau^*\phi$ are not congruent and that $\phi g_j$ and $(\tau^*\phi)g_j$ are not isometric follow almost exactly the proofs for $SO(2m+2r)$, by making the obvious adjustments for $Sp(2n+r)$.  Therefore our examples are non-trivial. 

\end{example}


\end{document}